\documentclass[12pt]{amsart}
\theoremstyle{definition}

\theoremstyle{plain}
\newtheorem{theorem}{Theorem}
\newtheorem{corollary}{Corollary}

\newtheorem{prop}{Proposition}
\newtheorem{ejem}{Example}

\theoremstyle{remark}

\DeclareMathOperator{\NN}{\mathbb N}

\title{Estimates of Disjoint Sequences in Banach Lattices and R.I. Function Spaces}

\author{R. Gonzalo and J.A. Jaramillo}

\address{Departamento de Matem\'atica Aplicada, Facultad de Inform\'atica de Madrid\\
         Universidad Polit\'ec\-ni\-ca, Campus de Montegancedo\\
         Boadilla del Monte, 28660 \\
         }

\email{rngonzalo@fi.upm.es}
\address{
         Departamento de An\'alisis Matem\'atico, Facultad de Ciencias Matem\'{a}ticas \\
         Universidad Complutense de Madrid. \\
28040 Madrid.
         }
\email{jaramil@mat.ucm.es}

\begin{document}

\begin{abstract}
We introduce $UDS_p$-property (resp. $UDT_q$-property) in Banach
lattices as the property that every  normalized disjoint sequence
has a subsequence with an upper $p$-estimate (resp. lower
$q$-estimate). In the case of rearrangement invariant spaces, the
relationships with Boyd indices of the space are studied. Some
applications of these properties are given to the high order
smoothness of Banach lattices, in the sense of the existence of
differentiable bump functions.


\end{abstract}

\maketitle

\section{Introduction}

Recall that, for $1 < p,q < \infty$, a sequence
$\{x_k\}_{k=1}^{\infty}$ in a Banach space $X$  is said to have an
{\it upper $p$-estimate} (respectively, a {\it lower
$q$-estimate}) with constant $C>0$ if for any $a_1,\dots, a_n\in
\Bbb R$ and $n \in \NN$,
$$
\Vert \sum_{k=1}^{n}a_k x_k \Vert  \leq C \left(\sum_{k=1}^{n}
\vert a_k \vert^p \right)^{1/p}
$$
$$
\left( \text{resp.}  \qquad C\left( \sum_{k=1}^{n} \vert a_k
\vert^q \right)^{1/q} \leq \Vert \sum_{k=1}^{n}a_k x_k \Vert
\right).
$$
According to \cite{K-O}, a Banach space $X$ is said to have {\it
property} $S_p$ if every weakly null normalized sequence in $X$
has a subsequence with an upper $p$-estimate. In a simi\-lar way,
$X$ is said to have {\it property} $T_q$  if every weakly null
normalized sequence has a subsequence with a lower $q$-estimate
(see \cite{G-JI}). These properties play an important role in the
behavior of polynomial maps on the space, as can be seen in
\cite{G-JI},\cite{CGG} and \cite{Dineen} and they are also related
to high order smoothness of the space ( see \cite{G-JI} and
\cite{G-JII}).

Here we study a variant of these properties in Banach lattices
and, more specifically, in the class of  rearrangement invariant
spaces. Note that in the simple case of spaces $L_p (\mu)$, where
$1<p<\infty$, since the space contains an isomorphic copy of
$\ell_2$, it verifies property $S_r$ just for $r \leq min\{ p,2
\}$. On the other hand, every {\it disjoint} normalized sequence
is equivalent to the unit vector basis of $\ell_p$. This remark
motivates the study of new properties related to the existence of
upper or lower estimates for {\it disjoint } sequences. For
technical reasons we are interested in the existence of such
estimates in an uniform way.

For $1 < p, q < \infty$ we say that a Banach lattice $X$ satisfies
the {\it uniform disjoint $S_p$-property}, in short
$UDS_p$-property, (respectively, {\it uniform disjoint
$T_q$-property, $UDT_q$- property} ) if there exists a constant
$C>0$ such that every normalized disjoint sequence in $X$
satisfies an upper $p$-estimate (resp. a lower $q$-estimate) with
constant $C$. It is then clear that spaces $L_p(\mu)$ for
$1<p<\infty$ verify both $UDS_p$ and $UDT_p$ properties.

To simplify the notation we introduce the following indices:
$$
\begin{array}{ll}
{\ell}_d(X)=& \sup \{ p> 1 : \quad X \quad  \text{has} \quad
UDS_p- \text{property}\} \\
 {\it u }_d(X)=& \inf \{ q >1 : \quad
X \quad \text{has} \quad UDT_q-\text{property}\}
\end{array}
$$

On the other hand, recall that for $1< p, q<\infty$, a Banach
lattice $X$ is said to have  a {\it lattice upper $p$-estimate}
(respectively, a {\it lattice lower $q$-estimate}) if there is a
constant $C>0$ such that, for every finite family of disjoint
vectors $x_1,\dots x_n$ in $X$, we have
$$ \Vert \sum_{k=1}^{n}
x_k \Vert \leq C \left( \sum_{k=1}^{n} \Vert x_k \Vert^p
\right)^{1/p}
$$
(respectively,
$$
\left( \sum_{k=1}^{n} \Vert x_k \Vert^q \right)^{1/q} \leq C\Vert
\sum_{k=1}^{n} x_k \Vert ).
$$
We also define:
$$
\begin{array}{ll}
up(X) = & \sup\{ p> 1 : \, X \, \text{has a lattice upper
p-estimate}\}
\\
low(X) =& \inf\{ q> 1 : \; X \, \text{has a lattice lower
q-estimate} \}.
\end{array}
$$

Of course, if $X$ has a lattice upper $p$-estimate (respectively,
a lattice lower $q$-estimate) then $X$ satisfies the
$UDS_p$-property (respectively,  $UDT_q$- property). As can be
seen with simple examples, the converse is not true. Indeed, for
$1<p,q< \infty$ consider the Banach lattice $X_{p,q} = \left(
\bigoplus_{n=1}^{\infty} \ell_q^n \right)_{\ell_p}$ with its
natural order. Using a standard sliding hump argument, it can be
shown that every normalized disjoint sequence in $X_{p,q}$ has a
subsequence equivalent to the unit vector basis of $\ell_p$, with
uniform constant. Then $X_{p,q}$ has both $UDS_p$ and $UDT_p$
properties. On the other hand, in the case $q<p$, $X_{p,q}$ has no
lattice upper $r$-estimate for $r>q$, and in the case $q>p$,
$X_{p,q}$ has no lattice lower $r$-estimate for $r<q$. Thus for a
Banach lattice $X$ we have that $up(X)  \leq \ell_d(X)$ and $
\it{u}_d(X) \leq low(X)$, and these inequalities are strict in
general.

In Section 2, we will relate these new properties $UDS_p$ and
$UDT_q$ with Boyd indices $p_X$ and $q_X$ in the context of
separable rearrangement invariant spaces $X(I)$, where $I$ is
either $\NN, [0,1]$ or $[0, \infty)$ (in short, r.i. spaces).
Recall that Boyd indices  play an important role in interpolation
problems, that is, the study of r.i. function spaces which are
between $L_p(I)$ and $L_q(I)$, in the sense that every operator
which is bounded in these two spaces, is bounded also on $X(I)$.
We prove that if $I=\Bbb N$ or $I=[0,\infty)$ and $X(I)$ has
$UDS_p$-property (respectively, $UDT_q$-property), then $p\leq
p_X$ (respectively, $q\geq q_X$). This fact is not true in the
case of r.i. function spaces on $[0,1]$ as Proposition $5$ below
will show.

Section 3 is devoted to study the connection between
$UDS_p$-property and high order smoothness of Banach lattices. The
smoothness of a space is here understood as the existence of {\it
bump} functions, that is, real valued functions with bounded
support, endowed with a certain degree of differentiability. In
this direction, in \cite{israel} it is proved that the existence
of a bump function which is $n$-times Fr\'{e}chet differentiable
with the last derivative $(p-n)$-H\"{o}lder continuous on a Banach
space $X$, not containing copies of $\ell_{2k}$, implies the
existence of upper $p$-estimates in the sequences of the space.
Again, since most of r.i. function spaces contain $\ell_2$ this
kind of results do not give very precise information on these
spaces. Our approach here is to consider  disjoint sequences in
the space. In this sense, we prove that if a Banach lattice $X$
does not contain  a vector sublattice order-isomorphic to
$\ell_{2k}$ and admits a  bump function which is $n$-times
Fr\'{e}chet differentiable with the last derivative
$(p-n)$-H\"{o}lder continuous, then $X$ has $UDS_p$-property. Some
related results are also obtained. Finally, we investigate the
best order of smoothness of Lorentz function spaces $L_{p,q}$.

\section{$UDS_p$, $UDT_p$-properties and Boyd indices.}

We first give some definitions. We refer to \cite{LT-II} for
details. Recall that if $X(I)$ is an r.i. function space on an
interval $I$ which is either $[0,1]$ or $[0,\infty)$, the Boyd
indices $p_X$, $q_X$ are defined by
$$
p_X= \lim_{s\to \infty} \frac{\log s}{\log \Vert D_s \Vert} \qquad
q_X= \lim_{s\to 0^+} \frac{\log s}{\log \Vert D_s \Vert},
$$
where for each $0<s<\infty$, the linear operator $D_s:X(I)\to
X(I)$ is defined by:

\noindent
\begin{enumerate}

\item In the case of $I=[0,\infty)$, for a measurable function $f$ on
$I$:
$$
(D_sf)(t)=f(t/s),\qquad  0\leq t<\infty
$$

\item In the case of $I=[0,1]$, for a measurable $f$ on $I$
:
$$
(D_sf)(t)= f(t/s) \chi_{[0,min\{1,s\}]}, ,\qquad  0\leq t \leq 1.
$$

\end{enumerate}

\noindent The Boyd indices can be also defined  for r.i. spaces on
the integers. In this case the operators $D_s$ are defined only if
$s$ is an integer or the reciprocal of an integer. If
$f=(a_1,a_2,\dots,)$ and $n=1,2,\dots$ we consider :
$$
D_nf=( a_1,\dots,a_1,a_2,\dots, a_2,\dots)
$$
$$
D_{1/n}f=n^{-1}\left( \sum_{i=1}^n a_i, \sum_{i=n+1}^{2n} a_i,
\dots \right).
$$

\noindent The indices $p_X$ and $q_X$ are defined as above by
taking the limits only over $s=n$  (respectively, $s=1/n$),
$n=1,2, \cdots$.

Now consider $X(\NN)$ a separable r.i.  space on $\NN$. Note that
in these spaces the basis $\{e_n\}=\{\chi_{\{n\}}\}$ is a
symmetric sequence, and therefore such spaces coincide with spaces
having symmetric basis. We have:

\begin{prop} Let $X(\NN)$ be a separable
r.i. space on $\NN$ and let $1<p,q < \infty$.
\begin{enumerate}
\item If $X(\NN)$ has $UDS_p$-property, then
$p\leq p_X$.
\item If $X(\NN)$ has $UDT_q$-property, then
$q\geq q_X$.
\end{enumerate}
\end{prop}

\begin{proof} We give the proof in the case of upper estimates; the
result for lower estimates is obtained in an analogous way.

\noindent Consider  $e_n=\chi_{\{n\}}$ for all $n\in \NN$. Since
the space has $UDS_p$-property, there is a constant $C>0$ such
that each  normalized disjoint sequence in $X(\NN)$ admits a
subsequence with an upper $p$-estimate with constant $C$. Let
$x=a_1e_1+\dots+a_ke_k\in X(\NN)$ with $\Vert x \Vert=1$. Then,
$$
D_n(x)=a_1(e_1+\dots +e_n)+\dots +a_k(e_{(n-1)k+1}+\dots +e_{nk})
= x_1+\dots + x_n
$$
where
$$
x_1 = a_1e_1+ a_2e_{n+1}+\dots +a_{k}e_{n(k-1)+1}
$$
  $$
x_n= a_1e_n+a_2e_{2n}+ \dots +a_ke_{kn}
$$
 and $\Vert x_1 \Vert = \dots = \Vert x_n \Vert=1$ .
\smallskip
Consider now, $x_{n+i}= a_1e_{(n+i-1)k+1}+ \dots +
a_ke_{(n+i-1)k+k}$ for $i=1,2,\dots$.
\smallskip
The sequence $\{x_n\}_{n=1}^{\infty}$ is normalized and disjoint
and therefore it admits a subsequence $\{x_{n_k}\}_{k=1}^{\infty}$
with an upper $p$-estimate with constant $C$. Thus,
$$
\Vert x_{k_1}+ \dots +x_{k_n} \Vert \leq Cn^{1/p}
$$
and consequently,
$$
\Vert x_1 + \dots +x_n \Vert = \Vert x_{k_1}+ \dots +x_{k_n} \Vert
\leq Cn^{1/p}
$$
Thus, $ \Vert D_n x \Vert \leq Cn^{1/p}$ where $C$ is independent
of $n$. Now, since the space is separable, sequences with finite
support are dense in the space, and therefore for every $n\in \NN$
we have $ \Vert D_n \Vert \leq Cn^{1/p}. $ It is easy to check now
that
$$
p_X = \lim_{n\to \infty} \frac{\log n}{\log \Vert D_n \Vert} \geq
p
$$
as we required.
\end{proof}

Now we would like to recall that, if a Banach lattice $X$ does not
contain a vector sublattice order-isomorphic to $\ell_1$, then
every normalized disjoint sequence in $X$ is weakly null. This is
a consequence (taking positive and negative parts) of the
following well-known result, which will be useful along the paper.
For a proof, see e.g \cite{AB}, Theorem 14.21.

\begin{prop} Let $\{x_n\}_{n=1}^{\infty}$ be a normalized disjoint sequence of positive vectors
in a Banach lattice $X$. The following are equivalent:
\begin{enumerate}
\item $\{x_n\}_{n=1}^{\infty}$ is not weakly null.
\item $\{x_n\}_{n=1}^{\infty}$ contains a subsequence equivalent to the unit vector basis of
$\ell_1$.
\end{enumerate}
\end{prop}

Suppose now that $X(\NN)$ is a separable r.i. space on $\NN$, not
containing any vector sublattice order-isomorphic to $\ell_1$.
Then every normalized disjoint sequence is weakly null and,
conversely, every weakly null normalized sequence is equivalent to
a disjoint sequence. On the other hand, by \cite{K-O}, if the
space has $S_p$-property then there is a constant $C>0$ such that
every weakly null sequence has a subsequence with an upper
$p$-estimate with the same constant $C$. Therefore in this case
$S_p$-property and $UDS_p$-property coincide, and $\ell_d(X)$
coincides with the index $\ell(X)= \sup \{ p \geq 1; \quad X \quad
\text{has} \quad S_p-\text{property} \} $ introduced in
\cite{G-JI}. That is, we have in this case $\ell(X)=\ell_d(X)\leq
p_X$. In some cases we have equality. For example, if $X=\ell_M$
is a separable Orlicz sequence space with $1<p_X<q_X<\infty$, it
is shown in \cite{raquel} that $\ell(\ell_M)=\alpha_X=p_X$ (see
also \cite{LT-II} Proposition 2.b.5). In general, $p_X$ is
strictly larger than $\ell(X)$, as the following example shows.

Recall that for $1<p,q< \infty$ the Lorentz sequence space
$X=\ell_{p,q}$ is defined as the space of all sequences
$\{a_n\}_{n=1}^{\infty}$ for which
$$
\Vert (a_{n})_{n=1}^{\infty}  \Vert_{p,q} = \left(
\sum_{n=1}^{\infty} (a_{n}^{*})^q(n^{q/p} - (n-1)^{q/p})
\right)^{1/q} < \infty
$$
\noindent where $\{a_{n}^{*}\}_{n=1}^{\infty}$ is the decreasing
rearrangement of $\{\vert a_n\vert\}_{n=1}^{\infty}$.

\begin{ejem} For $1< q<p<\infty$, the Lorentz sequence space $X= \ell_{p,q}$
satisfies $ \ell(X)=\ell_d(X)=q <p = p_X $.
\end{ejem}

\begin{proof} In the case $1< q<p<\infty$ we have that the Lorentz space
$\ell_{p,q}$  coincides with the Lorentz sequence space {\it d(w ;
q)} with weights $w_n=n^{q/p}-(n-1)^{q/p}$, and it is reflexive
(see \cite{LT-I} pg. 178). In particular, it does not contain any
copy of $\ell_1$. An easy computation shows that in these spaces
$\Vert D_n \Vert = n^{1/p}$ for each integer $n$. Consequently
both Boyd indices coincide with $p$, that is, $p_X=q_X=p$. On the
other hand, as shown in \cite{raquel}, $\ell(\it{d(w;q)})=q$.
\end{proof}

\medskip

Now we turn to the case of separable r.i. function spaces on an
infinite interval $X[0,\infty)$.

\begin{prop} Let $X[0,\infty)$ be a separable
r.i. function space on $[0,\infty)$ , and consider $1< p,q
<\infty$.

\begin{enumerate}
\item If $X[0,\infty)$ has $UDS_p$-property,
then $p\leq p_X$.
\item If $X[0,\infty)$ has $UDT_q$-property, then $q\geq  q_X$.
\end{enumerate}
\end{prop}
\begin{proof} We will prove (1), and (2) is analogous. Assume
that $X[0,\infty)$ has $UDS_p$-property with constant $C>0$. First
of all, since $X[0,\infty)$ is a separable r.i. function space,
step functions are dense on $X[0,\infty)$ and therefore, it is
enough to show that there is a constant $R>0$ such that
$$
\Vert D_s(f)  \Vert \leq Rs^{1/p} \qquad {\text if  } \quad s \geq
1
$$
for every norm-one step function $f$.

Now let $f=\sum_{i=1}^{m} \lambda_i \chi_{[u_{i-1}, u_{i})}$,
where $0=u_0<u_1<\dots <u_m=r$,  with $\Vert f \Vert=1$.  We
consider for all $k\in \NN$,
$$
f_n = \sum_{i=1}^m \lambda_i  \chi_{[(n-1)r+ u_{i-1}, (n-1)r+
u_i)}
$$
Then, the sequence $\{f_n\}_{n=1}^{\infty}$ constructed in such a
way is normalized and disjoint, and therefore it admits a
subsequence $\{f_{n_k}\}_{k=1}^{\infty}$ with an upper
$p$-estimate. Thus,
$$
\Vert f_{n_1} + \dots + f_{n_k} \Vert \leq Ck^{1/p}
$$
and
$$
\Vert f_1 +\dots +f_k \Vert = \Vert f_{n_1} + \dots + f_{n_k}
\Vert \leq Ck^{1/p}
$$
Moreover, with a suitable automorphism $\sigma$ of $[0,\infty)$ we
have that $f_{1} + \dots + f_k = g \circ \sigma $ where
$$
g= \lambda_1\chi_{[0,ku_1)} + \lambda_2\chi_{[ku_1,ku_2)}+ \dots +
\lambda_m \chi_{[ku_{m-1},ku_m)}.
$$
Now $g=D_kf$ since
$$
(D_kf)(t)= f(t/k).
$$
Thus,
$$
(D_kf)(kx)=f(kx/k)=f(x)
$$
and
$$
\Vert g \Vert = \Vert D_kf \Vert = \Vert f_1+ f_2+ \dots +f_k
\Vert \leq Ck^{1/p}
$$
for all $k\in \NN$. As a consequence,
$$
\Vert D_k f \Vert \leq Ck^{1/p}
$$
for all $k\in \NN$. Next, if $s>1$, there is an integer $k$ such
that $k \leq s <k+1$. Since $D_sf \leq D_{k+1}f $ we have that
$$
\Vert D_s f \Vert \leq \Vert D_{k+1}f \Vert \leq C(k+1)^{1/p}\Vert
f \Vert =
$$
$$
C(\frac{k+1}{k})^{1/p}k^{1/p}\Vert f \Vert \leq C2^{1/p}s^{1/p}
$$
Therefore, there is a constant $R>0$ such that
$$
\Vert D_s(f)  \Vert \leq Rs^{1/p}
$$
for any step function $f$ with $\Vert f \Vert =1$.
\smallskip

Now, it is easy to check that:
$$
p_X= \lim_{s\to \infty}\frac{\log s}{\log \Vert D_s \Vert} \geq p
$$
as we required.
\end{proof}

As a consequence, we obtain:

\begin{corollary} Let $X[0,\infty)$ be a separable r.i. function space. Then:

\begin{enumerate}
\item
 $up(X) \leq \ell_{d}(X) \leq p_X.$
\item
 $low(X) \geq {\it u}_d(X) \geq q_X.$
\end{enumerate}
\end{corollary}

For separable Orlicz function spaces $X=L_M[0,\infty)$ we obtain
equalities above. Indeed, it is shown in \cite{LT-II} (Proposition
2.b.5) that in this case $p_X=up(X)$ and $q_X=low(X)$.

On the other hand, in the case of separable Orlicz function spaces
$L_M[0,\infty)$, it is possible to obtain the following
characterization of $UDS_p$ and $UDT_q$ properties in terms of the
Orlicz function $M$:

\begin{prop} Let $L_M[0,\infty)$ be a separable Orlicz function space.
The following statements are equivalent:
\begin{enumerate}

\item $L_M[0,\infty)$ has $UDS_p$-property (respectively, $UDT_q$-property).

\item
$$
\inf_{t>0,\lambda \geq 1} \frac{M(\lambda t)}{\lambda^pM(t)} >0
\,\,\, \, \,  \, \, \, \, \, (resp. \sup_{t>0,\lambda \geq 1}
\frac{M(\lambda t)}{\lambda^qM(t)} < \infty)
$$
\end{enumerate}

\end{prop}

\begin{proof}
We prove the result concerning the $UDS_p$-property. The proof for
$UDT_q$-property is analogous. Proceeding as in the proof of
Proposition 3,  we have that there is a constant $R>0$ such that
$\Vert D_s \Vert \leq Rs^{1/p}$. In particular, if $u>0$ for
$f=\chi_{[0,u/s]}$ we have:
$$
\Vert \chi_{[0,u]} \Vert \leq Rs^{1/p}\Vert \chi_{[0,u/s]} \Vert
$$

We now proceed now as in \cite{LT-II} (pg. 139). Let $u>0$ and
$s\geq 1$ be fixed; since
$$
M^{-1}(\frac{1}{u})= \frac{1}{\Vert \chi_{[0,u]}\Vert}
$$

we have:
$$
M^{-1}\left(\frac{s}{u} \right) = \frac{1}{\Vert
\chi_{[0,u/s]}\Vert} \leq Rs^{1/p}\frac{1}{\chi_{[0,u]}}
=Rs^{1/p}M^{-1}(\frac{1}{u}).
$$

By taking: $t=M^{-1}(\frac{1}{u})$ and $\lambda=Rs^{1/p}$, it
follows that if $t>0$ and $\lambda \geq R$:
$$
M(\lambda t) \geq s/u = M(t)\lambda^p\frac{1}{R^{p}}
$$
and therefore,
$$
\inf_{t>0,\lambda \geq 1} \frac{M(\lambda t)}{M(t)\lambda^p} >0
$$
as we required.

In order to prove that (2) implies (1), by using \cite{LT-II} (see
pg. $140$) we obtain that if:
$$
\inf_{t>0,\lambda \geq 1} \frac{M(\lambda t)}{M(t) \lambda^p} >0
$$
then the space $L_{M}[0,\infty)$ has an upper $p$-estimate, and
consequently it has $UDS_p$-property.
\end{proof}

Now we give an application,  taking into account that, for
$1<p,q<\infty$ the function spaces $L_{p}[0,\infty) \cap
L_{q}[0,\infty)$ and $L_{p}[0,\infty)+L_q[0,\infty)$ can be
des\-cribed as Orlicz spaces. Namely, we have that
$L_{p}[0,\infty) \cap L_{q}[0,\infty)=L_M[0,\infty)$, where
$M(t)=\min \{t^p,t^q \}$ and
$L_{p}[0,\infty)+L_q[0,\infty)=L_M[0,\infty)$, where $M(t)=\max
\{t^p,t^q \}$. Therefore, with a simple calculation we obtain:

\begin{corollary} Let $1<p,q<\infty$. Then,
\begin{itemize}
\item[{\it(i)}] The spaces $L_{p}[0,\infty) \cap L_{q}[0,\infty)$ and
$L_{p}[0,\infty)+L_q[0,\infty)$ have $UDS_r$-property if and only
if $r \leq min\{p,q\}$.
\item[{\it(ii)}] The spaces $L_{p}[0,\infty) \cap L_{q}[0,\infty)$
and $L_{p}[0,\infty)+L_q[0,\infty)$ have $UDT_r$-property if and
only if $r\geq max\{p,q\}$.
\end{itemize}
\end{corollary}
\medskip

To finish this section, we see that the situation is quite
different in the case of r.i. function spaces on the interval
$[0,1]$. First recall that, if $1<p,q<\infty$ and $I$ is either
$[0,1]$ or $[0,\infty)$, the Lorentz space $L_{p,q}(I)$ is the
space of all measurable functions $f$ on $I$ for which
$$
\Vert f \Vert_{L_{p,q}} = \left( \int_I (f^{*}(t))^q
d(t^{\frac{q}{p}})\right)^{\frac{1}{q}} < \infty,
$$
\noindent where $f^*$ denotes the decreasing rearrangement of $f$.
We are grateful to E. Semenov for kindly providing us the
following result.

\begin{prop} Let $1<p,q<\infty$, and consider the space
$X=L_{p,q}[0,1]$. Then:

\begin{enumerate}

\item $X$ has Boyd indices $p_X = q_X = p$.

\item $X$ satisfies that ${\it up}(X) = \min\{p,q\}$ and ${\it low(X)} = \max\{p,q\}$.

\item $X$ has both $UDS_q$ and $UDT_q$-properties, and therefore
$$\ell_d(X)=u_d(X)=q.$$

\end{enumerate}
\end{prop}

\begin{proof} For $(1)$ we refer to \cite{BoI} and \cite{BoII}; and $(2)$
can be seen in \cite{Creekmore}. Now we pass to $(3)$. For the
case $1<q<p<\infty$, it is proved in \cite{Fiegel} (see Theorem
5.1) that, for each $\varepsilon
>0$, every disjoint normalized sequence in $L_{p,q}[0,1]$ admits a
subsequence which is $(1-\varepsilon)^{-1}$-equivalent to the unit
vector basis of $\ell_q$. Consequently, such space verifies
$UDS_q$ and $UDT_q$ properties. On the other hand, for the case
$1<q<p<\infty$ it is proved in  \cite{Carothers}  (see the proof
of Theorem 1),  that there is a constant $C>0$, only depending on
$p$ and $q$, in such a way that for each $\varepsilon>0$, every
normalized disjoint sequence in $L_{p,q}[0,1]$ admits an upper
$p$-estimate with constant $C=1-2\varepsilon$. Consequently, the
space has $UDS_q$-property. Moreover, the space admits a lattice
lower $q$-estimate with constant $q$. Consequently, $L_{p,q}[0,1]$
admits $UDS_q$ and $UDT_q$ properties.

\end{proof}

\section{Impact of high order of smoothness on upper estimates.}

Let $1<p<\infty$ and let $N$ be the greatest integer strictly less
than $p$. According to \cite{Meshkov}, a real valued function
$\phi:X\to \Bbb R$ is said to be $UH^p$-{\it smooth} if it is
$N$-times Fr\'{e}chet differentiable and there is a constant $C>0$
such that if $x,y \in X$ ,
$$
\Vert \phi^{(N)}(x) - \phi^{(N)}(y) \Vert \leq C \Vert
x-y\Vert^{p-N}.
$$
As usual, we say that the space $X$ is $UH^p$-smooth if there is a
$UH^p$-smooth  function on $X$ with bounded nonempty support.

In the case of separable r.i. spaces $X(\NN)$, the smoothness of
the space has a great impact on its structure as a Banach lattice.
More precisely, we have:

\begin{theorem}\cite{G-JII} Let $1<p< \infty $ and let $X$ be
a Banach space with a symmetric basis, not isomorphic to
$\ell_{2k}$ for $k \in \NN$. If $X$ is $UH^p$-smooth then it has a
lattice upper $p$-estimate.
\end{theorem}

This result does not hold in the general case of Banach lattices,
as the example $X=\left( \bigoplus_{n=1}^{\infty} \ell_4^n
\right)_{\ell_6}$ shows (see \cite{G-JII}, Remark 1.6). In order
to obtain a variant  of Theorem above for Banach lattices, we
study the relationship of high order smoothness of the space with
the property $UDS_p$ introduced in the previous section. We use
analogous techniques as those used in \cite{israel}. We say that a
sequence $\{x_n\}_{n=1}^{\infty}$ is $\mathcal{P}_N$-null for some
integer $N$, if $\lim_{n\to \infty} P(x_n)= 0$ for every
homogeneous polynomial on $X$ of degree at most $N$.

\begin{theorem} Let $1<p<\infty$
and let $X$ be a $UH^p$-smooth Banach lattice. Then, there is a
constant $C>0$ such that, if $\{x_n\}_{n=1}^{\infty}$ is a weakly
null, normalized and disjoint sequence in $X$, either
\begin{enumerate}
\item $\{x_n\}_{n=1}^{\infty}$ has a subsequence equivalent to the
basis of $\ell_{2k}$ with $2k\leq p$, or
\item $\{x_n\}_{n=1}^{\infty}$ has a subsequence with an upper
$p$-estimate with constant $C$.
\end{enumerate}
\end{theorem}

\begin{proof} Let $N$ be the greatest
integer strictly less than $p$. Let $\phi$ be a $UH^p$-smooth
function on $X$ verifying that $\phi(0)=0$ and $\phi(x)=2$ for
$\Vert x \Vert \geq 1$ (this is always possible composing with a
suitable smooth real function). Consider $M>0$ such that:
$$
\Vert \phi^{(N)}(x) - \phi^{(N)}(y) \Vert \leq M \Vert
x-y\Vert^{p-N}
$$
for $x,y \in X$.
\smallskip
Let $\{x_n\}_{n=1}^{\infty}$ be a weakly null normalized disjoint
sequence in $X$.  Two cases may be given:

\begin{enumerate}
\item The sequence $\{x_n\}_{n=1}^{\infty}$ is not $\mathcal{P}_N$-null.
Then let $r$ be the minimum of all $k \in \NN$ such that
$\{x_n\}_{n=1}^{\infty}$ is not $\mathcal{P}_k$-null, and note
that $1<r\leq N<p$. Since the sequence $\{x_n\}_{n=1}^{\infty}$ is
unconditional,  by Proposition 1.9 in \cite{israel} it admits a
subsequence with a lower $r$-estimate. On the other hand, by
Proposition 1.1. in \cite{israel}, there exists a subsequence with
an upper $r$-estimate, and consequently it is equivalent to the
unit vector basis of $\ell_{r}$. Since $X$ is $UH^p$-smooth and
contains $\ell_r$,  then $r$ must be an even integer $r=2k$ (see
e.g.  \cite{Bonic}). Therefore, the subsequence is equivalent to
the unit vector basis of  $\ell_{2k}$ with $2k\leq N$.

\item The sequence $\{x_n\}_{n=1}^{\infty}$ is
$\mathcal{P}_N$-null. Proceeding as in Proposition 1.1. in
\cite{israel} we obtain a subsequence $\{x_{n_i}\}_{i=1}^{\infty}$
verifying:
$$
\vert \phi(\sum_{i=1}^{k} a_if_{n_i}) - \phi(0) \vert \leq
M(\sum_{i=1}^{k} \vert a_i \vert^{p}) +1
$$
for all $a_1,\dots, a_k\in \Bbb R$ and $k \in \NN$. From this we
deduce that $\{x_n\}_{n=1}^{\infty}$ has a subsequence with an
upper $p$-estimate with constant $C=M^{1/p}$.
\end{enumerate}
\end{proof}

As a consequence we obtain the following:

\begin{corollary}
Let $1<p<\infty$ and let $X$ be a Banach lattice not containing
any vector sublattice order-isomorphic to $\ell_{2k}$ for $2k \leq
p$. If $X$ is $UH^p$-smooth then it has $UDS_p$-property.
\end{corollary}

\begin{proof}
Since $X$ is $UH^p$-smooth, it is super-reflexive (see e.g.
\cite{DGZ}, Theorem V.3.2), and in particular it contains no copy
of $\ell_1$. Now by Proposition $2$, every norma\-lized disjoint
sequence in $X$ is weakly null, and Theorem $2$ applies.
\end{proof}

Note that in the above corollary we cannot avoid the requirement
for $X$ of not containing a disjoint sequence equivalent to the
usual basis of $\ell_{2k}$. Indeed, for $X=L_{2k}(\mu)$, with
$k\in \NN$, we have that $X$ is $UH^p$-smooth for every $p>1$ but
$X$ has $UDS_p$-property only for $p\leq 2k$.

\begin{corollary} Let $X$ be a Banach lattice. If the index
$ \ell_{d}(X)$ is not an even integer, then the space is not
$UH^q$ -smooth for $q>\ell_{d}(X)$.
\end{corollary}

\begin{proof}
Assume the contrary. Suppose that  $\ell_d(X)$ is not  an even
integer and consider  $1<p<\ell_d(X)<q<\infty$ in such a way that
the interval $[p,q]$ does not contain any even integer, and $X$ is
$UH^q$-smooth. Let $\{x_n\}_{n=1}^{\infty}$ be a normalized
disjoint sequence in $X$. Since $X$ has $UDS_p$-property, there
exists a subsequence  with an upper $p$-estimate. Then,
$\{x_n\}_{n=1}^{\infty}$ has no subsequence equivalent to the
usual basis of $\ell_{2k}$ with $2k\leq p$. As in Corollary $3$,
since $X$ is $UH^q$-smooth it contains no copy of $\ell_1$, and by
Proposition $2$ we have that $\{x_n\}_{n=1}^{\infty}$ is weakly
null. Consequently, by Theorem $2$ there is a $C>0$ such that
$\{x_n\}_{n=1}^{\infty}$ has a subsequence with an upper
$q$-estimate with a constant $C$. Then, the space has
$UDS_q$-property, which is a contradiction.

\end{proof}

In \cite{Maleev} it is proved that for Orlicz spaces $\ell_M$,
$L_{M}[0,1]$ and $L_{M}[0,\infty)$, the smoothness properties of
the space depend on the Boyd indices. In the case of r.i. function
spaces on $[0,\infty)$, by using the results obtained in the
previous section, we obtain an estimate of the smoothness in terms
of the corresponding Boyd index:

\begin{corollary} Let $1<p<\infty$ and let $X[0,\infty)$
be a separable r.i. function space  not containing any vector
sublattice order-isomorphic to $\ell_{2k}$ for $2k\leq p$. If
$X[0,\infty)$ is $UH^p$-smooth then $p \leq p_X$.
\end{corollary}

In \cite{LT-II} it is proved that if $X$ is a r.i. function space
on $I=[0,\infty)$, then for every $1< p<p_X$ and $q_X<q<\infty$,
we have:
$$
L_p(I)\cap L_q(I) \subset X(I) \subset L_p(I)+L_q(I)
$$
with the inclusion maps being continuous. The converse of this
result is not true in general. As an application of Theorem $2$ we
may obtain the following result:

\begin{prop} Let $1<p<\infty$.
\begin{enumerate}
\item Suppose that $X(\NN) \subset \ell_p$ with continuous inclusion,
and $X(\NN)$ is not isomorphic to $\ell_{2k}$ for $2k\leq p$. Then
$X(\NN)$ is not $UH^r$-smooth for $r > p$.

\item Suppose that $X[0,\infty) \subset
L_{p}[0,\infty)+L_{q}[0,\infty)$ with continuous inclusion, and
$X[0,\infty)$ does not contain  any vector sublattice
order-isomorphic to $\ell_{2k}$ for $2k\leq \max\{p,q\}$. Then
$X[0,\infty)$ is not $UH^r$-smooth for $r>\max\{p,q\}$.
\end{enumerate}
\end{prop}

\begin{proof}
$(1)$ Note first that $X(\NN) \subset \ell_p$ with continuous
inclusion, means that the basis $\{e_n\}_{n=1}^{\infty}$ of
$X(\NN)$ has a lower $p$-estimate.  Assume that $X(\NN)$ is
$UH^r$-smooth with $p<r<[p]+1$. Since the basis
$\{e_n\}_{n=1}^{\infty}$ is equivalent to all of its subsequences,
by Theorem $2$ we obtain that $\{e_n\}_{n=1}^{\infty}$ has an
upper $r$-estimate, and this is not possible.

$(2)$ Suppose for example that $p<q$, and suppose that
$X[0,\infty)$ is $UH^r$-smooth for $r>q$. Consider the sequence
$f_n=\chi_{[n,n+1)}$ in $X[0,\infty)$. In the space
$L_{p}[0,\infty)+L_{q}[0,\infty)$, the sequence
$\{f_n\}_{n=1}^{\infty}$ is equivalent to the unit vector basis of
$\ell_q$. Then there exists some $C>0$ such that
$$
 C ( \sum_{n=1}^{\infty} \vert a_n \vert^q )^{1/q}  \leq \Vert
\sum_{n=1}^{\infty}  a_nf_n \Vert_{X}.
$$
This means that $\{f_n\}_{n=1}^{\infty}$ admits a lower
$q$-estimate and consequently $\{f_n\}_{n=1}^{\infty}$ does not
admit any upper $r$-estimate for $r>q$. This is a contradiction by
Theorem $2$.

\end{proof}

A similar result is the following.

\begin{prop}
Let $1<p<q<\infty$ and suppose that the interval $[p,q]$ does not
contain any even integer.
\begin{enumerate}

\item If $X(\NN)$ is a r.i. space and
$\ell_p \subset X(\NN) \subset \ell_q$ with continuous inclusions,
then $X(\NN)$ is not $UH^r$-smooth for $r>q$.

\item If $X[0,\infty)$ is a r.i. function space and
$L_{p}[0,\infty)\cap L_{q}[0,\infty) \subset X[0,\infty) \subset
L_{p}[0,\infty)+L_{q}[0,\infty)$ with continuous inclusions, then
$X[0,\infty)$ is not $UH^r$-smooth for $r>q$.

\end{enumerate}

\end{prop}
\begin{proof} $(1)$ Suppose that $X(\NN)$ is $UH^r$-smooth for $r>q$.
Note that since $\ell_p \subset X(\NN) \subset \ell_q$ with
continuous inclusions, we have that the symmetric basis
$\{e_n\}_{n=1}^{\infty}$ of $X(\NN)$, has both an upper
$p$-estimate and a lower $q$-estimate. Since $[p,q]$ does not
contain any even integer, $\{e_n\}_{n=1}^{\infty}$ has no
subsequence equivalent to the usual basis of $\ell_{2k}$ for any
$k\in \NN$. Therefore, we obtain a contradiction from Theorem $2$.

$(2)$  Suppose that $X[0,\infty)$ is $UH^r$-smooth for $r>q$.
Consider the sequence  $f_n= \chi_{[n,n+1)}$ in $X[0,\infty)$. In
the space $L_{p}[0,\infty)+L_{q}[0,\infty)$, the sequence
$\{f_n\}_{n=1}^{\infty}$ is equivalent to the unit vector basis of
$\ell_q$, and in the space $L_{p}[0,\infty)\cap L_{q}[0,\infty)$,
the sequence $\{f_n\}_{n=1}^{\infty}$ is equivalent to the unit
vector basis of $\ell_p$. Therefore there exist $C,D>0$ such that:
$$
C (\sum_{n=1}^{\infty} \vert a_n  \vert^{q})^{1/q} \leq \Vert
\sum_{n=1}^{\infty} a_n f_n \Vert_{X} \leq D ( \sum_{n=1}^{\infty}
a_n \vert^{p} )^{1/p} .
$$
That is, $\{f_n\}_{n=1}^{\infty}$  has both an upper $p$-estimate
and a lower $q$ estimate and, as before, we obtain a contradiction
with Theorem $2$.
\end{proof}

In the sequel we investigate the smoothness properties of the
Lorentz function spaces. The result for Lorentz sequence spaces
was already obtained in \cite{G-JII}.

\begin{prop} Let $1<p \neq q< \infty$.
\begin{enumerate}

\item The Lorentz sequence space $\ell_{p,q}$ and the Lorentz function
space $L_{p,q}[0,\infty)$ are not $UH^r$-smooth if $r
>\min\{p,q\}$.

\item The Lorentz space $L_{p,q}[0,1]$ is not $UH^r$-smooth if
$r>q$.
\end{enumerate}
\end{prop}

\begin{proof} Since $\ell_{p,q}$ is isometric to a
sublattice of $L_{p,q}[0,\infty)$, by \cite{G-JII} it follows that
these spaces are not $UH^r$-smooth if $r>\min\{p,q\}$.

For $L_{p,q}[0,1]$,  it follows  from  \cite{C-D} that every
disjoint normalized sequence has a subsequence which is equivalent
to the unit vector basis of $\ell_q$. In the case that $q$ is not
an even integer the result follows by using that $\ell_q$ is not
$UH^r$-smooth if $r>q$ (see e.g. \cite{Bonic}).

Otherwise, let $q$ be an even integer and assume that the space is
$UH^q$-smooth. We now prove that $L_{p,q}[0,1]$ is saturated with
subspaces of cotype $q$, that is, every infinite dimensional
subspace of $L_{p,q}[0,1]$ contains an infinite dimensional
subspace with cotype $q$. Indeed, let $X$ be an arbitrary infinite
dimensional subspace of $L_{p,q}[0,1]$ and choose
$\{f_n\}_{n=1}^{\infty}$ a weakly null normalized sequence in $X$,
which contains  no norm convergent subsequence. This is always
possible since $L_{p,q}[0,1]$ is super-reflexive. By using an idea
of Rabiger (see \cite{Rabiger}) one can easily prove that passing
to a subsequence $\{f_n\}_{n=1}^{\infty}$ either has a lower
$2$-estimate or it is equivalent to a disjoint sequence. In the
first case, proceeding as in \cite{Maleev} we obtain that the
sequence $\{f_n\}_{n=1}^{\infty}$ (passing to a subsequence if
necessary) has also an upper $2$-estimate; therefore such a
sequence is equivalent to the unit vector basis of $\ell_2$ and
generates a subspace with cotype $q$. In the second case, if the
sequence is disjoint, by \cite{C-D} it has a subsequence
equivalent to the unit vector basis of $\ell_q$. Since $\ell_q$
has cotype $q$ the result follows also in this case. Now,
according to \cite{Deville} the space $L_{p,q}[0,1]$ has a
separating polynomial. This is not possible, since by \cite{GGJ}
the only r.i. function spaces $X[0,1]$ which admit a separating
polynomial are the spaces $L_{2k}[0,1]$ up to an equivalent norm.
\end{proof}

\end{document}